\documentclass[a4paper, 11pt]{article} 

\usepackage[protrusion=true,expansion=true]{microtype} 
\usepackage{graphicx} 
\usepackage{wrapfig} 
\usepackage{changepage}
\usepackage{mathpazo} 
\usepackage[T1]{fontenc} 
\usepackage{amsmath}
\usepackage{amssymb}
\usepackage{graphicx,subfig}
\usepackage{caption}
\usepackage[numbers]{natbib}
\usepackage{listings}
\usepackage{color}

\definecolor{dkgreen}{rgb}{0,0.6,0}
\definecolor{gray}{rgb}{0.5,0.5,0.5}
\definecolor{mauve}{rgb}{0.58,0,0.82}

\makeatletter
\renewcommand\@biblabel[1]{\textbf{#1.}} 
\renewcommand{\@listI}{\itemsep=0pt} 

\renewcommand{\maketitle}{ 
\begin{flushright} 
{\LARGE\@title} 

\vspace{50pt} 

{\large\@author} 
\\\@date 

\vspace{40pt} 
\end{flushright}
}

\title{\textbf{An Introduction to Krylov Subspace Methods}\\ 
A Less Mathematical Way to Understand} 

\author{\textsc{Shitao Fan} 
\\{\textit{Zhejiang University}}} 

\date{\today} 

\begin{document}

\maketitle 



\begin{abstract}
  Nowadays, many fields of study are have to deal with large and sparse data matrixes, but the most important issue is finding the inverse of these matrixes. Thankfully,
  Krylov subspace methods can be used in solving these types of problem. However, it is difficult to understand mathematical principles behind these methods. In the first part of the article,
  Krylov methods are discussed in detail. Thus, readers equipped with a basic knowledge of linear algebra should be able to understand these methods. In this part, the knowledge of Krylov methods are put into some 
  examples for simple implementations of a commonly known Krylov method GMRES. In the second part, the article talks about CG iteration, a wildly known method which is very similar to Krylov methods. By comparison between CG iteration and Krylov methods, readers can get a better comprehension of Krylov methods based on CG iteration. In the third part of the article, aiming to improve the efficiency of Krylov methods, preconditioners are discussed. In addition, the restarting GMRES is briefly introduced to reduce the space consumption of Krylov methods in this part.
\end{abstract}

\hspace*{3,6mm}\textit{Keywords: Krylov subspace method, CG iteration, algorithm}  

\vspace{30pt} 

\newpage
\tableofcontents
\newpage
\section{Introduction}
When meeting a problem to solve a linear system $Ax=b$ where $A$ is a square matrix and $b$ is a vector, people might primarily consider the solution can be easy to be expressed by $A^{-1}b$. However, if $A$ is large,
calculation of the inverse of $A$ can be inefficient and space consuming. \\
In optimization, scholars usually think about the function of $A^{-1}$ to a vector, which means when inputting a vector $b$, by the function of $A^{-1}$, the system can calculate a $x$. There are many methods created to solve a linear system. However, if the matrix gets larger and sparser, some efficient ways should be found.\\
Here it comes to some methods connected with Krylov subspace. They are valid ways to deal with large and sparse matrixes. Krylov methods are more likely ways to reduce the dimension of the matrix. In other words, vulgar understanding of Krylov methods can be summarized by: \\
\textit{Inversion of large and sparse matrix $\Rightarrow$ Krylov method(usually GMRES) \\$\Rightarrow$Least linear square $\Rightarrow$ Inversion of small matrix}

\section{Krylov methods}
\subsection{Why use Krylov methods}
To solve a linear system, those familiar methods like Gauss elimination or other simple iteration(such as Jacobi iteration) may firstly put in to consideration. But confronted with colossal data matrix,
these methods are not as useful as dealing with middle and small size matrixes, because the consumption of Gauss elimination and Jacobi iteration is $O(n^{3})$.
Moreover, these matrixes might be indefinite and asymmetrical. In this situation, using Krylov methods is really effective to conquer a large and linear system.

\subsection{What is Krylov subspace}
Krylov-based algorithms iteratively compute an approximation of $x$ by transforming an $n-$dimensional vector space into a lower $m-$dimensional($m\leq n$) subspace using matrix-vector multiplications without requiring estimation of $A^{-1}$ explicitly. Usually, get a rough estimate of $Ax=b$ by $x_{0}$, then let $x = x_{0}+x_{m}$. The estimate  of $x_{m}$ forms an approximation in the Krylov subspace $\mathcal{K}_{m}$:\\
$$\mathcal{K}_{m}(A,r_{0}) = span\{r_{0},Ar_{0},...,A^{m}r_{0}\},\quad r_{0} = b-Ax_{0}$$
$"Span"$ means every vector in $\mathcal{K}_{m}(A,r_{0})$ can be expressed by the linear combination of the basis $  \{r_{0},Ar_{0},...,A^{m}r_{0}\}$.\\
In other words, a conclusion can be drew: $x_{m} \thickapprox \sum_{i=1}^{m-1}\beta_{i}A^{i}r_{0}$.\\
If eager to estimate $x_{m}$ accurately, a $k$ must be found to meet: \\
$$\mathop {Az }\limits_{z \in \mathcal{K}_{k}(A,r_{0})} = \mathop{y}\limits_{y \in \mathcal{K}_{k}(A,r_{0})}$$
$\mathcal{K}_{k}(A,r_{0}) = span\{r_{0},Ar_{0},...,A^{m}r_{0}\}$ is called the invariant subspace of $A$. However, $k$ is connected with the scale of the matrix. When the matrix gets larger, calculation becomes more onerous.
Thus, a small number $m$ are usually set to make the error $r_{m}= b-A(x_{0}+x_{m})$ smaller than the tolerance.
\subsection{How to use Krylov methods}
\subsubsection{Algorithm of Krylov methods}
This section concentrates on Arnoldi's Method at first, which is a way to generalize a vector group to span Krylov subspace. to begin with, find a normalized vector in Krylov subspace. Then, use $A$ to multiply the
vector and exploit Schmidt orthogonalization making the vector perpendicular to the formal one. After m times of that procedure, Arnoldi's Method generalizes m vectors which are perpendicular to each other.\\

\textbf{ALGORITHM Arnoldi}\\
1. \qquad choose a vector $\nu_{1}$, $s.t.$ $||\nu_{1}||_{2}=1$.\\
2. \qquad For $j = 1,2,..,m$ Do:\\
3. \qquad \qquad Compute $h_{ij} = (A\nu_{j},\nu_{i})$ for $i = 1,2,..,j$\\
4. \qquad \qquad Compute $\omega_{j} = A\upsilon_{j}-\sum_{i=1}^jh_{ij}\nu_{i}$\\
5. \qquad \qquad $h_{j+1,j} = ||\omega_{j}||_{2}$\\
6. \qquad \qquad $if h_{j+1,j} = 0 then stop$\\
7. \qquad EndDO\\
\\
In this algorithm, Arnoldi method generalizes a unit orthogonal basis $V_{m} = (\nu_{1},...,\nu_{m})$ of Krylov subspace. Plus, noticing the matrix $H_{m}$ consisting of $h_{ij}$. It is not hard to draw a conclusion
(The proof of it is written at appendix Proof1) :
$$AV_{m} = V_{m}H_{m}$$
Then, return to solve the linear system $Ax=b$. When using Krylov method with less calculation, errors are inevitable. Therefore, Krylov methods want to find a small $m$ to get the estimate of $x$ by $x_{0}+x_{m}$
and make the residue $r_{m}$ as small as possible. It is described in the following equation:
$$\mathop {\min }\limits_{x_{m}\in \mathcal{K}_{m}(A,r_{0})}||Ax-b||_{2}, x=x_{0}+x_{m}$$
Let $r_{0} = b-Ax_{0}, \beta:= ||r_{0}||_{2},$ , $\nu_{1}:=r_{0}/\beta$ and $e_{1}=[1,0,0,...,0]^{T}$\\
Estimate $||Ax-b||_{2}$ as \\
$||Ax-b||_{2} \\
= ||A(x_{0}+x_{m})-b||_{2} \\
= ||A(x_{0}+V_{m}y)-b||_{2}\\
= ||AV_{m}y-r_{0}||_{2}\\
=||V_{m}H_{m}y-r_{0}||_{2}\\
=||V_{m}H_{m}y-\beta \nu_{1}||_{2}\\
=||V_{m}({H}_{m}y-\beta e_{1}) ||_{2}$\\
So the former problem converts into:
$$\mathop {\min }\limits_{y}||{H}_{m}y-\beta e_{1}||_{2}$$
And in this way, the reduction dimension property of Krylov methods indicated in the introduction section has been implemented.  And if $m\ll $ the scale of $A$, these steps have already converted the problem to be a
simple one with less calculation. The following algorithm is a way to implement it.\\
\\
\textbf{ALGORITHM GMRES}\\
1.\qquad Compute $r_{0} = b-Ax_{0}, \beta:= ||r_{0}||_{2},$ and $v_{1}:=r_{0}/\beta$\\
\textit{get a rough estimate of the equation and initialize the prim vector of krylov subspace.}\\
2.\qquad For $j = 1,2,..,m$ Do:\\
3.\qquad \qquad Compute $\omega_{j} = A\upsilon_{j}$ \\
\textit{construct the next vector to span krylov subspace.}\\
4.\qquad \qquad For $i = 1, . . ., j$ Do:\\
5.\qquad \qquad \qquad $h_{ij} = (\omega_{j},\upsilon_{i})$\\
6.\qquad \qquad \qquad $\omega_{j} = h_{ij}\upsilon_{j}$\\
\textit{Schmidt orthogonize the vector.}\\
7.\qquad \qquad EndDo\\
8.\qquad \qquad $h_{j+1,j} = ||\omega_{j}||_{2}$, $if h_{j+1,j} = 0$, set $m = j$ and break\\
\textit{test whether the subspace is invariant suspace.}\\
9.\qquad \qquad $\upsilon_{j+1} = \omega_{j}/h_{j+1,j}$\\
\textit{normalize the vector}\\
10.\qquad EndDo\\
11.\qquad Define the $(m + 1)\ast m$ Hessenberg matrix $H_{m} = \{h_{ij}\}_{ 1\leq i\leq m+1,1\leq j\leq m}$\\
12.\qquad Compute $y_{m}$ the minimizer of $||\beta e_{1} - {H}_{m}y_{m}||_{2}$ and $x_{m} = x_{0} + V_{m}y_{m}$.\\
\textit{use least linear square to minimize the residue.}
\subsubsection{Simple implement of GMRES}
In this section, a $3\times3$ matrix is used as an example to show how GMRES works. Its matlab code is written in the appendix name Code1.\\
\\
EXAMPLE1 GMRES\\
$$
 A=\left[
 \begin{matrix}
   1 & 4 & 7 \\
   2 & 9 & 7 \\
   5 & 8 & 3
  \end{matrix}
  \right]
  b=\left[
 \begin{matrix}
   1  \\
   8  \\
   2
  \end{matrix}
  \right]
$$
Let $x_{0} = 0$, $r_{0} = b$, and $maxit = 3$\\
$V_{m} = (v_{1},v_{2},v_{3})$\\
$H_{m} = \{h_{i,j}\}_{(m+1)\times m} $\\
$\beta:= ||r_{0}||_{2} = 8.31,\quad \nu_{1} = r_{0}/\beta = [0.12,0.96,0.24]^{T}$\\
$\omega_{1}=A\nu_{1} = [5.66,10.59,9.03]^{T}$\\
$h_{1,1} = (\omega_{1},\nu_{1})= 13.06$\\
$\omega_{1} = \omega_{1}-h_{1,1} \nu_{1} = [4.09,-198,5.88]^{T}$\\
$h_{2,1} = ||\omega_{1}||_{2} = 7.43$\\
$\upsilon_{2} = \omega_{1}/h_{2,1}= [0.55,-0.27,0.79]^{T}$\\
In the same way, all the vectors and matrixes can be calculated:\\
$$
 V_{m}=\left[
 \begin{matrix}
   0.12 & 0.55 & 0.82 \\
   0.96 & -0.27 & 0.037 \\
   0.24 & 0.79 & -0.56
  \end{matrix}
  \right]
  H_{m}=\left[
 \begin{matrix}
   13 & 5.4& -1.6  \\
   7.4 & 4.0 &1.1  \\
   0 & 2.6 &-4.1 \\
   0 & 0 & 0
  \end{matrix}
  \right]
$$
So now the problem changes to $argmin||\beta e_{1} - {H}_{m}y_{m}||_{2}$. It is easy to get:
$$y_{m} = {[1.36,-2.16,-1.4]}^{T}, x = 0 + x_{m} = {[-2.18, 1.84, -0.6]}^{T}$$
\newpage
\noindent
EXAMPLE2 GMRES\\
Like what EXAMPLE1 does, EXAMPLE3 uses the same wat solve a linear system: $AX=B$
$$
 A=\left[
 \begin{matrix}
   1 & 4 & 7 \\
   2 & 9 & 7 \\
   5 & 8 & 3
  \end{matrix}
  \right]
  b=\left[
 \begin{matrix}
   1 & 2 & 5 \\
   8 & 3 & -3 \\
   2 & 9 & 8
  \end{matrix}
  \right]
$$
To solve it, $X$ should be separated by $(x_{1},x_{2},x_{3})$ and $B$ should be separated by $(b_{1},b_{2},b_{3})$. Thus the problem is divided into $Ax_{i}=b_{i}, \quad i=1,2,3$And the result of it is:
$$
 X=\left[
 \begin{matrix}
   -2.2 & 2.1 & 4.8 \\
   1.8 & -0.22 & -2.6 \\
   -0.59 & 0.11 & 1.5
  \end{matrix}
  \right]
$$
\\
EXAMPLE3 GMRES\\
In this example, the reduction dimension property will be talked here.\\
Here A is a random $100 \times 100$ matrix, and b is a random $100 \times 1$ vector. Use GMRES function in matlab, then carry out the sentences:
\begin{lstlisting}[title=Use GMRES to solve full matrix, frame=shadowbox]
A =  randi([10 100],100,100);
if det(A)~=0
    b = randi([10 100],100,1);
    maxit=90;
    tol=1e-1;
    [x2,fl2,rr2,it2,rv2] = gmres(A,b,[],tol,maxit);
end
\end{lstlisting}
The residue of $Ax=b$ is still bigger than $0.1$, which indicates that GMRES is not proper to solve full matrixes. In addition, the reduction dimension property of GMRES is hard to realize based on a full matrix.\\
\\
However if use GMRES to deal with a sparse matrix $A$, whose size is $7585 \times 7585$ matrix, and b is a $7585 \times 1$ vector. Use the same program, things will be different: \\
The solution of GMRES converges to $10^{-7}$ after 164 times iterations, and the residue is less than $9.8 \times 10^{-8}$. Here, the reduction dimension property of GMRES is showed by converting the $7585 \times 7585$ matrix to $164 \times 164$ scale of problem.\\
\section{Comparison between CG(conjugate gradient) method and Krylov subspace method-GMRES}
\subsection{Conjugate gradient method}
This section will briefly introduce a widely known method CG iteration, and compare it with GMRES. It is meaningful to focus on CG method. Not only it is efficient way to deal with positive definite symmetric and sparse linear system, but it is easily to be used.\\
At first, CG method is put forward to solve the optimization problem $min$ $f(x)=\frac{1}{2}x^{T}Ax-bx^{T}$, which is equal to find the answer of $g(x) = f^{'}(x) = Ax-b = 0$\\
To abridge the algorithm, the iteration sequence is named as $\{x_{k}\}$ and $x_{k+1} = x_{k} + \alpha_{k} d_{k} $. Also use the signal $g_{k}$ to represent $g(x_{k})$ or $f^{'}(x_{k})$.\\
\textbf{ALGORITHM CG}\\
1.\qquad find a rough estimate of the equation, such as $x_{0}$, and let $k=0$\\
\textit{initialize the iteration}\\
2.\qquad $d_{0} = g_{0}$\\
\textit{find the fastest downward direction as the first direction}\\
3.\qquad For $k=1,2,...,n-1:$\\
4.\qquad \qquad if $g_{k} =0$ break\\
\textit{test whether get the exact answer}\\
5.$\qquad \qquad\alpha_{k}=\frac{g_{k}^{T}d_{k}}{d_{k}^{T}Gd_{k}}$\\
\textit{search the step length exactly}\\
6.$\qquad \qquad x_{k+1} = x_{k} + \alpha_{k} d_{k}$\\
7.$\qquad \qquad g_{k+1} = Ax_{k+1}-b$\\
8.$\qquad \qquad\beta_{k+1} = \frac{g_{k+1}^{T}Ad_{k}}{d_{k}^{T}Gd_{k}}$\\
9.$\qquad \qquad d_{k+1}=-g_{k+1}+\beta_{k+1}d_{k}$\\
\textit{calculate the next conjugate direction}\\
\subsection{Comparison}
\textbf{Similarity:} To begin with, both of CG and GMRES are finding a estimate of $x_{m}$ where $Ax_{m} = b-r_{m}-Ax_{0}$. By using projection to subspace or subspace spanned by conjugate vector group, these two methods evaluate $x_{m}$ by the linear combination of the basis of subspace (like figure 1) .\\
So giving an $n \times n$ positive definite matrix, GMRES and CG can get the accurate answer by n-step iteration.\\
\begin{figure}[h!]
\centering
\includegraphics[width = 8cm,height=8cm]{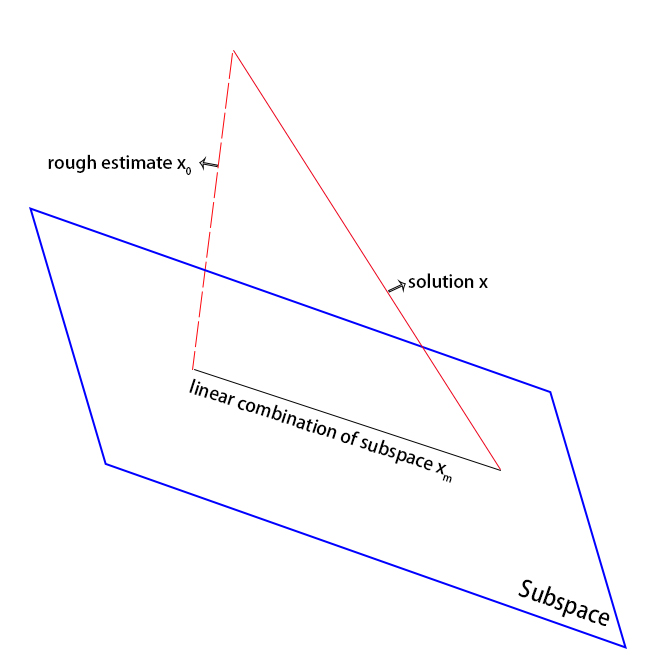}\\
\caption{projection}
\label{1}
\end{figure}
\\
\textbf{Difference:} CG iteration can only deal with symmetrical and positive definite matrix, but GMRES can easily solve asymmetrical matrix.\\
\\
Here let us know some of new concepts for further understanding:\\
Use $p(A)$ to represent $\sum_{i=1}^{k-1}\beta_{i}A^{i}$.\\
$$\mathcal{K}(A) = \frac{\lambda_{max}}{\lambda_{min}}$$
The value of $\mathcal{K}(A)$ is connected with $cond(A)$, and if the $cond(A)$ is large meaning the matrix A is ill.\\
The convergence property of CG and GMRES are different. The following part talks about when and how the two methods converge:\\
(Let $x^{*}$ to be the accurate solution of $Ax=b$ , $q(A)=1-tp(A)$, $p(A)$ is the linear combination of the basis in subspace, $||x^{0}-x^{*}||_{A} = \epsilon _{0}^{T}{q(A)}^{T}Aq(A)\epsilon _{0}$. Let $\lambda$ to be the eigenvalues of matrix $A$, $\Lambda$ to be the eigenvalue set of matrix $A$, $\mathbb{P}$ to be the polynomial ring)\\
CG iteration: $\frac{||x^{m}-x^{*}||_{A}}{||x^{0}-x^{*}||_{A}}\leqslant 2 (\frac{\sqrt{\mathcal{K}(A)}-1}{\sqrt{\mathcal{K}(A)}+1})^{2}$ \\
It means if $\mathcal{K}(A)\rightarrow 0$ the CG iteration will converge slowly, but if $\mathcal{K}(A)\rightarrow 1$ it will converge quickly.\\
GMRES: $\frac{||b-Ax^{m}||_{2}}{||r_{0}||_{2}}\leq min max|q(\lambda)|$, $q \in \mathbb{P}$ and $\lambda \in \Lambda$\\
The convergence of GMRES is really connected with eigenvalues of matrix $A$. Therefore it is not so easy to judge whether a matrix has a good property of convergence. However, there is some conclusion
that this method is superlinear. Thus, in normal sense GMRES will convergence faster than CG.\\                                                                                                                          \section{Further introduce about GMRES}
\subsection{Precpnditioner applied to GMRES and CG}
\subsubsection{What is preconditioner}
Firstly, the linear system is: $Ax=b$. However, the condition of A is bad, maybe ill or too large to calculate. So something should be found to improve the condition. The problem are changed by multiply a
preconditioner into:\\
\centerline{${M}^{-1}Ax = {M}^{-1}b$  or  $A{M}^{-1}u = {M}^{-1}b, x={M}^{-1}u$}\\
They have the same answer of the linear system, but the matrix changes to ${M}^{-1}A$ or $A{M}^{-1}$. The bad condition can be ameliorated.\\
Also, a common situation is when the preconditioner is available in the factored form,
$$M = M_{L}M_{R}$$
where typically $M_{L}$ and $M_{R}$ are triangular matrices. In this situation, the preconditioning can be split:
$${M_{L}}^{-1}A{M_{R}}^{-1}u = {M_{L}}^{-1}b, x={M_{R}}^{-1}u$$
\subsubsection{Why use preconditioner in GMRES and CG}
Preconditioner can help to improve the $cond$ of a matrix $A$ by converting to $cond(M^{-1}A)$. Also, it can really reduce the calculation of a large matrix.
That is, it will make the matrix be suitable for the method be suitable for the algorithm.\\
\subsubsection{How to find preconditioner in GMRES and CG}
For GMRES, faced with the asymmetrical sparse matrix. A method, incomplete LU factorization (ilu), is often put into consideration firstly. That is, find an upper triangle matrix $L$ and a lower triangle matrix $M$.
$A$ is separated by:\\
\centerline{$A=LU-R$, $R$ is the residue}\\
Assuming the number of the nonzero elements of the matrix is $n_{A}$, the consumption of the ilu is approximate $O({n_{A}}^\frac{3}{2})$. Compared with Gauss elimination, ilu does not need lots of calculation.\\
\\
For CG iteration, confronted with symmetrical sparse matrix is common. Incomplete Cholesky factorization (icho) is widely used as a preconditioner. In that process, separate $A$ by $L^{T}L$ where $L $ is a
lower triangle matrix.\\
\centerline{$A=LL^{T}-R$, $R$ is the residue}\\
Conspicuously, icho needs less calculation than ilu. Thus, it is also a valid way to find a preconditioner.
\subsubsection{How to use preconditoner in GMRES and CG}
For example, just consider the left-preconditioned method. The problem now is ${M}^{-1}Ax = {M}^{-1}b$.\\
When using GMRES and CG, simply consider ${M}^{-1}A$ as $A$ and ${M}^{-1}b$ as $b$. The algorithm is the same as GMRES.\\
\\
EXAMPLE4 GMRES AND CG WITH PRECONDITIONER\\ In this example, a large matrix is used here to demonstrate the effect of preconditioner.\\
The matrix here is:
\begin{figure}[h!]
\centering
\includegraphics[width = 10cm]{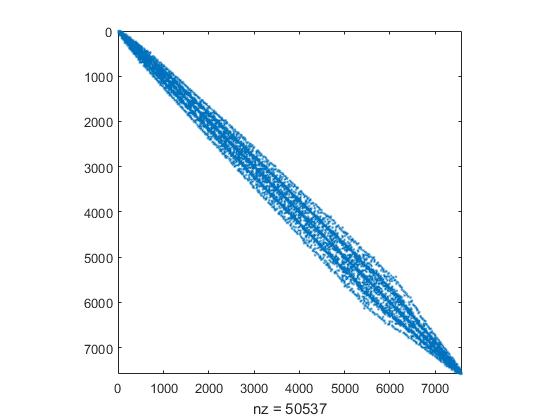}
\caption{The sparse matrix}
\label{2}
\end{figure}
\newpage
\noindent 
By use a special ilu method, IC(0), to find the preconditioner of matrix $A$. Two figures is displayed in the following part to show the effect of preconditioner on CG and GMRES. 
\begin{figure}[h!]
\centering
\includegraphics[height=8cm,width = 12cm]{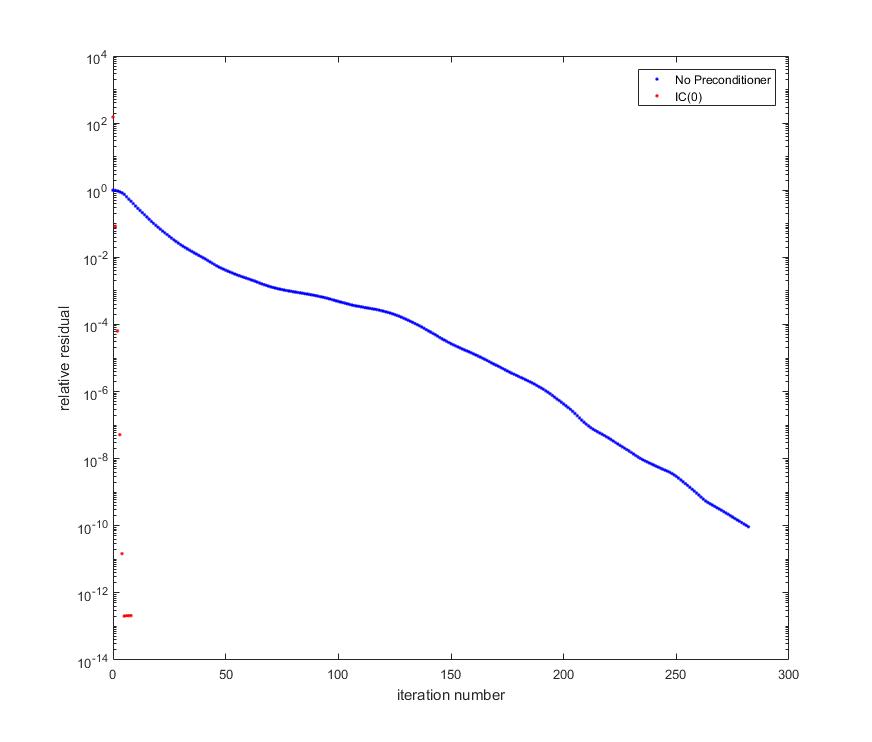}
\caption{GMRES with preconditioner}
\label{3}
\end{figure}
\begin{figure}[h!]
\centering
\includegraphics[height=8cm,width = 12cm]{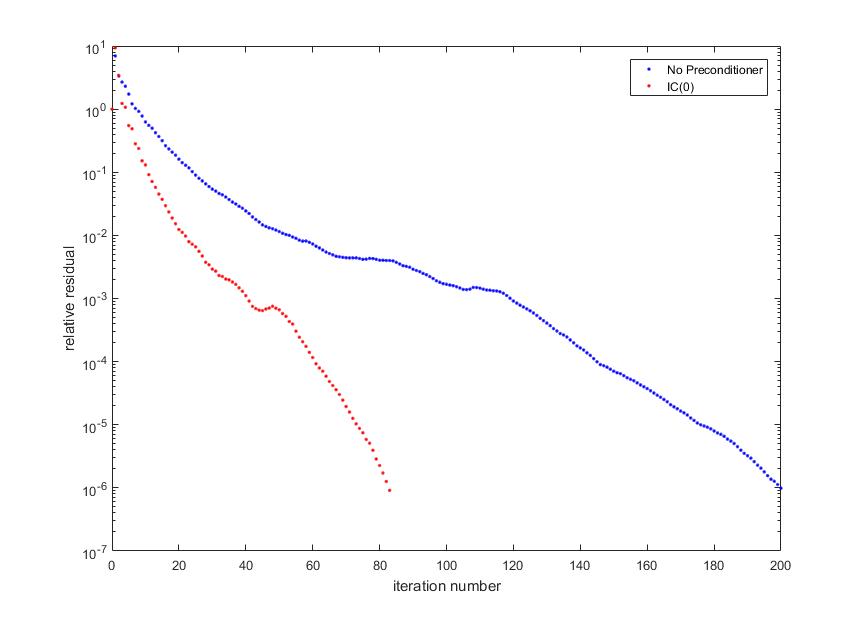}
\caption{CG with preconditioner}
\label{4}
\end{figure}
\newpage
\subsection{Restarting of GMRES}
If eager to get a more accurate estimate $\hat{x}$ of $Ax=b$, a larger m should be set. When the scale of $A$ is large enough, though $A$ is sparse, the vector group $V_{m}=(\nu_{1},...,\nu_{m})$ is not sparse. Thus,
the storage of $V_{m}$ is space consuming. A restart of GMRES can easily solve the problem. \\
That is, find a small $m_{1}$, and use GMRES get the answer $x_{m_{1}}$ after $m_{1}$ times of iteration. Then use $x_{m_{1}}$ as the iteration initial value. Repeat it until it converges to the tolerance. At each
step, just  a small matrix $V_{m_{1}}$ should be stored.\\                                                                                                                                                               \section{Appendix}
\subsection{Proof1}
\textit{proof:} By the process of Arnoldi¡¯s Method, a equation can easily be drew:
$$A\nu_{j} = \sum^{j+1}_{i=1}h_{ij}\nu_{i}, \quad j=1,2,...,m$$
Write $V_{m}$ by $\{\nu_{1},..,\nu_{m}\}$, then:
\begin{figure}[h!]
\centering
\includegraphics[width = 8cm]{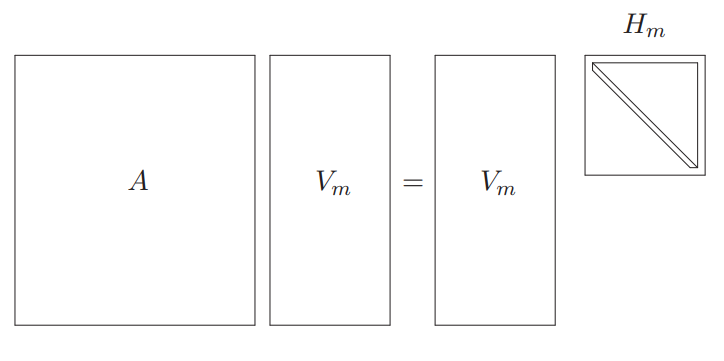}
\caption{Matrix multiplication}
\label{5}
\end{figure}
\newpage
\subsection{Code1}
\begin{lstlisting}[title=GMRES, frame=shadowbox]
function [xs,ys,Vs,Hs] =gmres_my(A,b,maxit)
r0=b;
beta=norm(r0);
V(:,1)=r0/norm(r0);
H=[];
for j=1:maxit
    W(:,j)=A*V(:,j);
    for i=1:j
        H(i,j)=dot(W(:,j),V(:,i));
        W(:,j)=W(:,j)-H(i,j)*V(:,i);
    end
    H(j+1,j)=norm(W(:,j));
    if H(j+1,j)==0
        break;
    end
    V(:,j+1)=W(:,j)/H(j+1,j);
end
[n,m]=size(H);
a=zeros(n,1);
a(1)=beta;
ys=lsqlin(H,a)
for i=1:length(ys)

    xs(i)=V(i,1:length(ys))*ys;
end
Vs=V;
Hs=H;
\end{lstlisting}


\bibliographystyle{unsrt}

\bibliography{sample}


\end{document}